\documentclass[12pt]{amsart}
\usepackage{mathrsfs}
\usepackage{amsfonts}
\usepackage{amssymb}
\usepackage{amsfonts, amscd, amsmath, mathrsfs, amssymb, amsthm, amsxtra, bbding, epsfig, eucal, graphicx, latexsym, url, mathbbol, bbold}
\usepackage[papersize={7.6in,10.2in},textwidth=14.7cm,textheight=21cm,centering]{geometry}
\usepackage{enumerate}

\usepackage[colorlinks=true,citecolor=blue,linkcolor=blue]{hyperref}
\hypersetup{
pdfstartpage=1,
pdfstartview=FitH}

\DeclareFontFamily{U}{mathx}{\hyphenchar\font45}
\DeclareFontShape{U}{mathx}{m}{n}{
      <5> <6> <7> <8> <9> <10>
      <10.95> <12> <14.4> <17.28> <20.74> <24.88>
      mathx10
      }{}
\DeclareSymbolFont{mathx}{U}{mathx}{m}{n}
\DeclareMathAccent{\widecheck}{\mathalpha}{mathx}{"71}

 \usepackage{caption} 
\numberwithin{equation}{section}

\allowdisplaybreaks

\newtheorem{lemma}{Lemma}[section]

\makeatletter
\newcounter{roem}
\renewcommand{\theroem}{\Roman{roem}}

\newcommand{\c@org@eq}{}
\let\c@org@eq\c@equation
\newcommand{\org@theeq}{}
\let\org@theeq\theequation

\newcommand{\setroem}{
\let\c@equation\c@roem
 \let\theequation\theroem}

\newcommand{\setarab}{
\let\c@equation\c@org@eq
\let\theequation\org@theeq}
\makeatother

\newtheorem*{claim*}{Claim}

\theoremstyle{remark}

\newcommand{\ue}{\mathrm{e}}

\DeclareMathOperator{\Mod}{mod}

\renewcommand{\bmod}[1]{\,(\Mod{ #1})}

\newcommand{\ba}{\boldsymbol{a}}

\newcommand{\bn}{\boldsymbol{n}}

\newcommand{\bZ}{\mathbf{Z}}

\newcommand{\cB}{\mathcal{B}}

\newcommand{\cR}{\mathcal{R}}

\newcommand{\tc}{\tilde{c}}
\newcommand{\tm}{\tilde{m}}
\newcommand{\tn}{\tilde{n}}

\def\le{\leqslant}

\def\ge{\geqslant}

\usepackage{graphicx}
\usepackage{tikz}

\begin{document}

\vglue -2mm

\title[Addendum to $q$-vdC]{Addendum to ``Arithmetic exponent pairs of algebraic trace functions and applications"}
\author{Jie Wu}

\address{%
School of Mathematics and Statistics
\\
Qingdao University
\\
308 Ningxia Road
\\
Qingdao
\\
Shandong 266071
\\
China}
\curraddr{%
CNRS, UMR 8050\\
Laboratoire d'Analyse et de Math\'ematiques Appliqu\'ees\\
Universit\'e Paris-Est Cr\'eteil\\
61 Avenue du G\'en\'eral de Gaulle\\
94010 Cr\'eteil cedex\\
France
}
\email{jie.wu@math.cnrs.fr}

\author{Ping Xi}

\address{School of Mathematics and Statistics, Xi'an Jiaotong University, Xi'an 710049, China}

\email{ping.xi@xjtu.edu.cn}

\maketitle

This addendum devotes to a detailed proof for the inequality (9.14) in our joint work \cite{WX21}.
We do not intend to publish this addendum in any journals; arXiv should be a good place for those reader who want to find such details.
The proof involves various averages of the following arithmetic functions.

\begin{itemize}

\item For $n\geqslant1$, denote by $n^\flat$ and $n^\sharp$ the squarefree and squarefull parts of $n$, respectively;
 i.e.,
\begin{align*}
n^\flat= \prod_{p\| n} p,\ \ \  n^\sharp= \prod_{p^{\nu}\| n, \, \nu\ge 2} p^{\nu}.
\end{align*}

\item 
For $n\geqslant1$, define
\begin{align*}
n^\ddagger=\prod_{p^2\|n}p.
\end{align*}

\item 
For $n\geqslant1$, define
\begin{align*}
\Xi(n)= \prod_{p^{\nu}\| n, \, \nu\ge 3} p^{\nu}.
\end{align*}
\end{itemize}

We keep the same convention as in \cite{WX21}, and the citation numbers are also referred to \cite{WX21} directly. 

\section{Preliminary lemmas}

For each positive integer $n$, we may write
\[n=n^\flat\cdot n^\sharp
\qquad\text{and}\qquad
(n^\flat, n^\sharp)=1.\]
The following lemma shows that $n^\sharp$ and $\Xi(n)$
are bounded by $\log n$ in certain averaged senses, and we will use this frequently without any notifications.

\begin{lemma}\label{lm:squarefull}
Let $x\geqslant2.$

{\rm (a)}
For any fixed $t\leqslant\frac{1}{2},$ we have
\begin{align*}
\sum_{n\leqslant x}(n^\sharp)^t\ll x\log x
\qquad\text{and}\qquad
\sum_{n\leqslant x}(n^\flat)^{-t}\ll x^{1-t}\log x.
\end{align*}
\par
{\rm (b)}
For any fixed $t\leqslant\frac{3}{4},$ we have
\begin{align*}
\sum_{n\leqslant x}\Xi(n)^t\ll x\log x.
\end{align*}
\end{lemma}
\proof 
Since $n\mapsto n^\sharp$ is multiplicative,
for $\Re e\, s>1,$ we have
\begin{align*}
\sum_{n\geqslant1} (n^\sharp)^t n^{-s}
& = \prod_p (1+p^{-s}+p^{-2(s-t)}+p^{-3(s-t)}+\cdots)
\\
& = \zeta(s)\zeta(2s-2t)A_t(s),
\end{align*}
where the Euler product $A_t(s)$ converges absolutely for $\Re e\,s>\max\{t+\frac{1}{3},\frac{1}{2}\}$
and $A_t(s)\ll 1$ for $\Re e\,s\ge \frac{6}{7}$ and $t\le \tfrac{1}{2}$. 
The first statement of assertion (a) follows immediately from a routine application of Perron's formula.
The second inequality is a consequence of the first one. 

As before, for $\Re e\,s>1,$ we have
\begin{align*}
\sum_{n\geqslant1} \Xi(n)^t n^{-s}
= \zeta(s)\zeta(3s-3t)B_t(s),
\end{align*}
where $B_t(s)$ is holomorphic and bounded for $\Re e\,s\ge1-\frac{1}{2}\varepsilon,t\le \tfrac{4}{5}-\varepsilon$.
A simple application of Perron's formula leads to the required inequality.
\endproof

We also need the following estimate on the greatest common divisors on average.

\begin{lemma}\label{lm:gcd}
For each $q\geqslant1$ and all $x\geqslant1$, we have
\begin{align*}
\sum_{n\leqslant x}(n,q)\leqslant\tau(q)x.
\end{align*}
\end{lemma}

\proof 
The sum in question can be rewritten as
\begin{align*}
\sum_{d|q}d\sum_{\substack{n\leqslant x\\(n,q)=d}}1
=\sum_{d|q}d\sum_{\substack{n\leqslant x/d\\(n,q/d)=1}} 1
\le \sum_{d|q,d\leqslant x}d\cdot x/d\le \tau(q)x,
\end{align*}
as claimed.
\endproof

The following inequality allows one to separate variables in the function $\Xi(\cdot).$
\begin{lemma}\label{lm:ddagger}
Let $m$ be squarefree and $s\geqslant 1$.
Then
\begin{align*}
\Xi(ms)\leqslant (s^\ddagger)^3 \Xi(s)^{\frac{4}{3}}.
\end{align*}
\end{lemma}

\proof
We can write uniquely $s=s_1 s_2^2 s_3$, 
where $s_1$ is squarefree, $s_2$ is squarefree and $s_3$ is cubefull such that $s_1, s_2, s_3$ are pairwise coprime.
We may also write uniquely $m=m_1m_2m_3m_4$ with $m_i\mid s_i^{\infty}$ for $1\le i\le 3$ and $(m_4, s)=1$. Thus we have
$$
\Xi(ms)
= \Xi(m_2s_2^2\cdot m_3s_3)
= \Xi(m_2s_2^2) \Xi(m_3s_3)
\leqslant (s^\ddagger)^3 \Xi(s_3)^{\frac{4}{3}}\leqslant (s^\ddagger)^3 \Xi(s)^{\frac{4}{3}}
$$
as claimed.
\endproof

\section{Initial treatment for $\cB(M,N)$}

Suppose $m$ and $n$ are squarefree numbers such that $(m, n)=1$.
Their divisors are also squarefree and pairwise coprime.
In what follows, we keep these in mind without mentioning specially for the sake of simplicity.
Moreover, all convention and notation in this section are all directly borrowed from \cite{WX21}.

Put
$$\Delta:=h_2\sigma_2a_1-h_1\sigma_1a_2,\quad
k := (\Delta, (\tm\tn s)^\flat).
$$
We would like to state two alternative estimates for 
\begin{align*}
T(\ba; m, \bn, s) 
= \sideset{}{^*}\sum_{r\sim R_0} 
\mathrm{e}(\xi([a_1,a_2]r)),
\end{align*}
in two complementary ranges, where $\xi$ is a fraction given by (9.10).

(I) If $R\geqslant q\delta[a_1,a_2]=\tm\tn s[a_1, a_2]/k,$ the completing method gives
\begin{align}\label{eq:T-estimate1}
T(\ba; m, \bn, s)
\ll X^\varepsilon R_0(\Delta, (\tm\tn s)^{\ddagger})\sqrt{\frac{\Xi(\tm\tn s) k}{\tm\tn s}}.
\end{align}

(II) If $R<q\delta[a_1,a_2]=\tm\tn s[a_1, a_2]/k$,
\begin{align}\label{eq:T-estimatebyexponentpairs}
T(\ba; m, \bn, s)
\ll X^{\varepsilon}
\bigg(\frac{(\tm\tn s)^\flat/k}{R_0}\bigg)^{\kappa} R_0^\lambda \big((\tm\tn s)^\sharp\big)^\nu,
\end{align}
where $(\kappa,\lambda,\nu)$ is an exponent pair defined by $(\boldsymbol{\varOmega})$ subject to an additional restriction 
corresponding to (7.1) or (7.2).

Together with the upper bound \eqref{eq:T-estimatebyexponentpairs},
the restrictions (7.1) and (7.2) produce inequalities involving $(\tm\tn,s^\flat),$
$(\Delta,(\tm\tn s)^\flat)$ and so on. To simply our arguments, we would like to introduce a trick on enlarging moduli. More precisely, suppose we are given an exponential sum
\[\sum_{n\in I}\ue\Big(\frac{f_1(n)\overline{f_2(n)}}{Q}\Big),\]
where $f_1,f_2\in\bZ[X]$ and zeros of $f_2\bmod Q$ have been excluded from summation. We aim to give a non-trivial upper bound for the above sum in terms of $Q$ and $|I|$. In practice, $Q$ might be not of a good shape, but $mQ$ can look much better for some positive integer $m$. Equivalently, we turn to consider the average
\[\sum_{n\in I}\ue\Big(\frac{\widetilde{f_1}(n)\overline{f_2(n)}}{mQ}\Big),\]
where $\widetilde{f_1}=mf_1.$
This is of course a trivial observation. While estimating exponential sums, it is a common feature that the problem becomes more difficult if the length of summation becomes shorter, or equivalently, the modulus becomes larger. The introduction of $m$ pushes the problem to a more difficult one at a first glance, but we need it to eliminate the bad structure of $Q$ and in the meanwhile, $m$ can also be controlled effectively (at least in an averaged sense).

Back to our applications to estimating $T(\ba; m, \bn, s),$ recall that the original modulus is $(\tm\tn s)^\flat/k$. We now enlarge the modulus to $q_{\text{new}}:=D\cdot(\tm\tn s)^\flat/k$ with $D=b(n_1c_2,n_2c_1)(\tm\tn,s^\flat)k/(\Delta,s^\flat).$ Precisely speaking, we may write
\begin{align*}
T(\ba; m, \bn, s) 
&= \sideset{}{^*}\sum_{r\sim R_0} \ue\Big(\frac{r\Delta[a_1,a_2]}{k}\frac{\overline{a_1a_2(\tm\tn s)^\sharp(r^2[a_1,a_2]^2+s^2)}}{(\tm\tn s)^\flat/k}\Big)W_\delta(r)\\
&= \sideset{}{^*}\sum_{r\sim R_0} \ue\Big(\frac{rD\Delta[a_1,a_2]}{k}\frac{\overline{a_1a_2(\tm\tn s)^\sharp(r^2[a_1,a_2]^2+s^2)}}{q_{\text{new}}}\Big)W_\delta(r).
\end{align*}
In this way, the method of arithmetic exponent pairs is in fact applied to the trace function
\begin{align}\label{eq:tracefunction-enlargedmoduli}
x\mapsto 
\mathrm{e}\bigg(\frac{xD\Delta[a_1,a_2]}{k}\frac{\overline{a_1a_2(\tm\tn s)^\sharp([a_1,a_2]^2x^2+s^2)}}{p}\bigg).
\end{align}

Note that 
\[\frac{k}{(\Delta,s^\flat)}=\frac{(\Delta,(\tm\tn s)^\flat)}{(\Delta,s^\flat)}=\frac{(\Delta,\tm\tn s^\flat/(\tm\tn,s^\flat))}{(\Delta,s^\flat)}=(\Delta,\tm\tn/(\tm\tn,s^\flat)),\]
from which one may see that $k/(\Delta,s^\flat)$ is a divisor of $\tm\tn.$ The latter one comes from the well-factorable weight in the linear sieve of Iwaniec, so that it consists of very small prime factors, which we mean {\it allowing good factorizations}.
 Similarly, 
the factor $b(n_1c_2,n_2c_1)(\tm\tn,s^\flat)$ also allows good factorizations.
On the other hand,
\begin{align*}
q_{\text{new}}=\frac{b(n_1c_2,n_2c_1)(\tm\tn,s^\flat)(\tm\tn s)^\flat}{(\Delta,s^\flat)}
=b(n_1c_2,n_2c_1)\tm\tn \cdot \frac{s^\flat}{(\Delta,s^\flat)},
\end{align*}
and one may find only the factor $s^\flat/(\Delta,s^\flat)$ does not allow good factorizations. 

We are now in a good position to applying Theorem 7.1 to estimating $T(\ba; m, \bn, s)$ with the trace function given by \eqref{eq:tracefunction-enlargedmoduli} and the deformation factor $W_\delta(\cdot)$, getting
\begin{align}\label{eq:T-estimate2}
T(\ba; m, \bn, s)
\ll X^{\varepsilon}
\bigg(\frac{\tm\tn s^\flat b(n_1c_2,n_2c_1)}{R_0(\Delta,s^\flat)}\bigg)^{\kappa} R_0^\lambda \big((\tm\tn s)^\sharp\big)^\nu,
\end{align}
for which the restrictions corresponding to (7.1) and (7.2) become 
\begin{align*}\tag*{$[\cR_1]$}
\Big(b(n_1c_2,n_2c_1)\tm\tn \cdot \frac{s^\flat}{(\Delta,s^\flat)}\Big)^{1-2\kappa}R_0^{2\kappa-2\lambda+1}\geqslant \frac{s^\flat}{(\Delta,s^\flat)}
\end{align*}
and
\begin{align*}\tag*{$[\cR_2]$}
\Big(b(n_1c_2,n_2c_1)\tm\tn \cdot \frac{s^\flat}{(\Delta,s^\flat)}\Big)^{2-2\lambda}R_0^{2\lambda-2\kappa-1}\geqslant \frac{s^\flat}{(\Delta,s^\flat)},
\end{align*}
respectively. As we will see later, the enlargement of the moduli will allow us to transform the above restrictions more freely.

In fact, when $\Delta=0,$ one can bound $T(\ba; m, \bn, s)$ trivially, so that the relevant contribution to $\cB(M,N)$
is 
\begin{align*}
& \ll S\sum_{\substack{n_1\sim N\\ n_2\sim N}} 
\sum_{\substack{h_1\le H\\ h_2\le H}}\sum_{m\sim M}\sum_{\substack{a_1\mid mn_1\\a_2\mid mn_2}}\sum_{\substack{r\sim R/[a_1,a_2]\\h_2a_1n_1/(n_1,r)=h_1a_2n_2/(n_2,r)}}1\\
&\ll S\sum_{n_1\sim N} 
\sum_{h_2\le H}\sum_{m\sim M}\sum_{a_1\mid mn_1}\sum_{r\sim R}\tau_3(h_2a_1n_1/(n_1,r))\\
&\ll HMNRSX^\varepsilon\\
&\ll HLMNX^\varepsilon.
\end{align*}
Henceforth, we always assume $\Delta\neq0$, in which case
denote by $\cB_1$ and $\cB_2$ the relevant contributions to $\cB(M,N)$ deriving from the above two estimates \eqref{eq:T-estimate1} and \eqref{eq:T-estimate2}, respectively. Hence we have
\begin{align}\label{eq:B(M,N)-decomposition}
\cB(M,N)=\cB_1+\cB_2+O(HLMNX^\varepsilon+H^3N).
\end{align}

\section{Estimate for $\cB_1$}
Throughout this section, we keep in mind the restriction
\begin{align}\label{eq:restriction1}
R<\tm\tn s[a_1, a_2]/k
\end{align}
in the summand.
Since $[a_1, a_2] = [bc_1, bc_2] = b[c_1, c_2]$, we can write
\begin{align*}
R_0(\Delta, (\tm\tn s)^{\ddagger})\sqrt{\frac{\Xi(\tm\tn s) k}{\tm\tn s}}
& = R \frac{(\Delta, (\tm\tn s)^{\ddagger})}{b[c_1, c_2]}
\sqrt{\frac{(\Delta, (\tm\tn s)^{\flat}) \Xi(\tm\tn s)}{\tm\tn s}}.
\end{align*}
Thus
$$
\cB_1
\ll X^\varepsilon R \sum_{b\le 2M}
\sum_{\substack{c_1\le 2N\\c_2\le 2N}} 
\sum_{\substack{h_1\le H\\ h_2\le H}} 
\sum_{\substack{m\sim M/b\\ s\sim S}} 
\sum_{\substack{n_1\sim N/c_1\\ n_2\sim N/c_2}} \frac{(\Delta, (\tm\tn s)^{\ddagger})}{b[c_1, c_2]}
\sqrt{\frac{(\Delta, (\tm\tn s)^{\flat}) \Xi(\tm\tn s)}{\tm\tn s}},
$$
where 
$$
\tn:=\frac{n_1n_2}{(n_1, n_2)}=[n_1, n_2],
\qquad
\Delta 
= \frac{b(h_2c_1n_1-h_1c_2n_2)}{(n_1, n_2)}.
$$
Furthermore, we write 
\begin{align*}
&c=(c_1, c_2),
\qquad
c_1=c\tc_1,
\qquad
c_2=c\tc_2,
\qquad 
(\tc_1, \tc_2)=1,
\\
& n=(n_1, n_2),
\qquad
n_1=n\tn_1,
\qquad
n_2=n\tn_2,
\qquad 
(\tn_1, \tn_2)=1.
\end{align*}
Then
\begin{align*}
\cB_1
&\ll X^\varepsilon R 
\sum_{\substack{b\le 2M\\ c\le 2N}} 
\sum_{\substack{\tc_1\le 2N/\tilde{c}\\ \tc_2\le 2N/c}} 
\sum_{\substack{h_1\le H\\ h_2\le H}} 
\sum_{\substack{m\sim M/b\\ n\le \min\{2N/c\tc_1,2N/c\tc_2\}\\ s\sim S}}\\
&\ \ \ \ \times
\sum_{\substack{\tn_1\sim N/nc\tc_1\\ \tn_2\sim N/nc\tc_2}} 
\frac{(\Delta, (m\tn s)^{\ddagger}) }{bc\tc_1\tc_2} 
\sqrt{\frac{(\Delta, (mn\tn_1\tn_2 s)^{\flat}) \Xi(mn\tn_1\tn_2 s)}{mn\tn_1\tn_2 s}},
\end{align*}
where $\Delta$ becomes $\Delta=bc(\tc_1h_2\tn_1-\tc_2h_1\tn_2)$. Note that the restriction \eqref{eq:restriction1} becomes
\begin{align*}
(\Delta,(mn\tn_1\tn_2 s)^\flat)\geqslant\frac{mn\tn_1\tn_2 sbc\tc_1\tc_2}{R}.
\end{align*}

Using the trivial inequality $q^j\le \sum_{d\mid q} d^j$ 
and Lemma \ref{lm:ddagger} with $(m, s)\mapsto(mn\tn_1\tn_2, s)$, 
we may have
\begin{align*}
\cB_1
&\ll X^\varepsilon R 
\sum_{\substack{b\le 2M\\ c\le 2N}} 
\sum_{\substack{\tc_1\le 2N/c\\ \tc_2\le 2N/c}} 
\sum_{\substack{m\sim M/b\\ n\le \min\{2N/c\tc_1,2N/c\tc_2\}\\ s\sim S}} 
\sum_{\substack{\tn_1\sim N/nc\tc_1\\ \tn_2\sim N/nc\tc_2}} 
\sum_{\substack{d_1\mid (mn\tn_1\tn_2 s)^{\flat}\\ d_2\mid (mn\tn_1\tn_2 s)^{\ddagger}\\d_1>mn\tn_1\tn_2 sbc\tc_1\tc_2/R}} \\
&\ \ \ \ \ \times\sum_{\substack{h_1\le H, \, h_2\le H\\ d_1d_2\mid bc(\tc_1h_2\tn_1-\tc_2h_1\tn_2)\neq0}} 
 \frac{d_1^{\frac{1}{2}} d_2 (s^\ddagger)^{\frac{3}{2}} \Xi(s)^{\frac{2}{3}}}{bc \tc_1\tc_2\sqrt{mn\tn_1\tn_2 s}}.
\end{align*}
After a dyadic partition, it follows that, with $N_0:=\min\{2N/c C_1, 2N/c C_2\}$, 
\begin{equation}\label{UB:B22_A}
\begin{aligned}
\cB_1
& \ll X^\varepsilon N^{-1}R 
\sum_{\substack{b\le 2M\\ c\le 2N}} 
\sum_{\substack{C_1\\ C_2}} 
\sum_{\substack{\tc_1\sim C_1\\ \tc_2\sim C_2}} 
\sum_{\substack{m\sim M/b\\ n\le N_0\\ s\sim S}} 
\sum_{\substack{\tn_1\sim N/nc C_1\\ \tn_2\sim N/nc C_2}} 
\\\noalign{\vskip -1mm}
&\ \ \ \ \ \times
\sum_{\substack{d_1\mid (mn\tn_1\tn_2 s)^{\flat}\\ d_2\mid (mn\tn_1\tn_2 s)^{\ddagger}\\d_1>mn\tn_1\tn_2 sbc\tc_1\tc_2/R}} 
\sum_{\substack{h_1\le H, \, h_2\le H\\ d_1d_2\mid bc(\tc_1h_2\tn_1-\tc_2h_1\tn_2)\not=0}} 
\frac{\sqrt{nd_1} d_2 (s^\ddagger)^{\frac{3}{2}} \Xi(s)^{\frac{2}{3}}}{b\sqrt{C_1C_2ms}},
\end{aligned}
\end{equation}
where we have used the fact that $\tn = n\tn_1\tn_2\asymp N^2/nc^2C_1C_2$ 
and $C_1,C_2$ are both of the shape $2^j$ with $0\le j\le\log(2N).$

Due to the symmetry between $(\tc_1,h_2,\tn_1)$ and $(\tc_2,h_1,\tn_2)$, we may impose the restriction
$\tc_1h_2\tn_1>\tc_2h_1\tn_2$ without loss of generality. For simplicity, the relevant contribution will also be denoted by $\cB_1$.

The next step is to sum over $h_2,n_1$ with necessary congruence conditions. To do so, we relax the conditions $d_1\mid (mn\tn_1\tn_2 s)^\flat$ and $d_2\mid (mn\tn_1\tn_2 s)^{\ddagger}$ firstly.

\vskip 1mm

\noindent{$\bullet$ \textit{Condition} $d_1\mid (m\tn s)^{\flat}$}

\vskip 1mm

Note that 
\[
(mn\tn_1\tn_2 s)^{\flat} 
= mn\tn_1\tn_2 \frac{s^\flat}{(mn\tn_1\tn_2, s^\flat)}\cdot\]
We may decomposite $d_1$ by $d_1=d_{1,1}d_{1,2}d_{1,3}d_{1,4}d_{1,5}$
with
\begin{align*}
d_{1,1}\mid m,
\qquad 
d_{1,2}\mid n,
\qquad 
d_{1,3}\mid \tn_1,
\qquad 
d_{1,4}\mid \tn_2,
\qquad 
d_{1,5}\Big|\frac{s^\flat}{(mn\tn_1\tn_2, s^\flat)}\cdot
\end{align*}
Since $d_1d_2\mid bc(\tc_1h_2\tn_1-\tc_2h_1\tn_2)$ and $(bc\tn_2, d_{1,3})=1$, 
we must have $d_{1,3}\mid \tc_2h_1$.
Thus, we may replace the condition of summation $d_1\mid (mn\tn_1\tn_2 s)^\flat$ in \eqref{UB:B22_A}
by $d_{1,1}$, $d_{1,2}$, $d_{1,3}$, $d_{1,4}$, $d_{1,5}$ with the restrictions
$$
d_{1,1}\mid m,
\qquad 
d_{1,2}\mid n,
\qquad 
d_{1,3}\mid \tc_2h_1,
\qquad 
d_{1,4}\mid \tn_2,
\qquad 
d_{1,5}\mid s^\flat.
$$

\vskip 1mm

\noindent{$\bullet$ \textit{Condition} $d_2\mid (m\tn s)^{\ddagger}$}

\vskip 1mm

Note that 
\[
(mn\tn_1\tn_2 s)^{\ddagger}\mid (m, s^\flat) (n, s^\flat) (\tn_1, s^\flat) (\tn_2, s^\flat){s^{\ddagger}}.
\]
We also have $(d_2,bc)=1$ since all prime factors of $d_2$ come from $s$ and $(bc,s)=1.$ 
Therefore, the congruence condition $d_1d_2\mid bc(\tc_1h_2\tn_1-\tc_2h_1\tn_2)$
is equivalent to 
$$
\tc_1h_2\tn_1\equiv \tc_2h_1\tn_2 \bmod{d_1d_2/(d_1, bc)}.
$$

In what follows, we will utilize the original notation $d_1,d_2$ and keep in mind the above decompositions. In this way, we can derive that
\begin{equation}\label{UB:B22_B}
\cB_1
\ll X^\varepsilon N^{-1} R 
\sum_{\substack{b\le 2M\\ c\le 2N}} 
\sum_{\substack{C_1\\ C_2}} 
\sum_{\substack{m\sim M/b\\ n\le N_0\\ s\sim S}} 
\sum_{\substack{\tc_2\sim C_2\\ h_1\le H\\ \tn_2\sim N/nc C_2}}
\sum_{\substack{(d_1)\\ (d_2)}} \frac{\sqrt{nd_1} d_2 (s^\ddagger)^{\frac{3}{2}} \Xi(s)^{\frac{2}{3}}}{b\sqrt{C_1C_2ms}} 
\mathcal{H},
\end{equation}
where
$$
\mathcal{H}
:= \sum_{\substack{\tc_1\sim C_1,\, h_2\le H, \, \tn_1\sim N/nc C_1\\ 
d_1d_2\mid bc(\tc_1h_2\tn_1-\tc_2h_1\tn_2)>0}} 1.
$$
Since $\tc_1h_2\tn_1>\tc_2h_1\tn_2$, the sum $\mathcal{H}$ is empty unless 
$d_1d_2/(d_1, b\tc)\ll HN/(\tc n),$
in which case one has
\begin{align*}
\mathcal{H}
= \sum_{\substack{\tc_2h_1\tn_2<d\le HN/nc\\ d\equiv \tc_2h_1\tn_2 \bmod{d_1d_2/(d_1, b\tc)}}}
\sum_{\substack{\tc_1\sim C_1, \, h_2\le H, \, \tn_1\sim N/nc C_1\\ 
\tc_1h_2\tn_1=d}} 1
& \ll X^\varepsilon \frac{HN}{c d_1d_2 n}(d_1, bc),
\end{align*}
from which and \eqref{UB:B22_B}, and noting the previous restriction $d_1>mn\tn_1\tn_2 sbc\tc_1\tc_2/R,$ we find
$$
\cB_1
\ll X^\varepsilon HR 
\sum_{\substack{b\le 2M\\ c\le 2N}} 
\sum_{\substack{C_1\\ C_2}} 
\sum_{\substack{m\sim M/b\\ n\le N_0\\ s\sim S}} 
\sum_{\substack{\tc_2\sim C_2\\ h_1\le H\\ \tn_2\sim N/nc C_2}}
\frac{(s^\ddagger)^{\frac{3}{2}} \Xi(s)^{\frac{2}{3}}}{bc \sqrt{C_1C_2mns}} \Big(1+\frac{MN^2S}{nc R}\Big)^{-\frac{1}{2}}
\sum_{(d_1)}(d_1, bc).
$$

Since $(d_{1,2}d_{1,3}d_{1,4}d_{1,5}, bc)=1$, we have
$(d_1, bc)=(d_{1,1}, bc).$
On the other hand, we trivially have
\begin{align*}
\sum_{d_{1,1}\mid m}(d_{1,1}, bc)\leqslant \tau(m)(m,bc).
\end{align*}
Summing over $m,\tc_2,h_1,\tn_2$ firstly, it follows that
\begin{align}
\cB_1
&\ll X^\varepsilon HR 
\sum_{\substack{b\le 2M\\ c\le 2N}} 
\sum_{\substack{C_1\\ C_2}} 
\sum_{\substack{m\sim M/b\\ n\le N_0\\ s\sim S}} 
\sum_{\substack{\tc_2\sim C_2\\ h_1\le H\\ \tn_2\sim N/nc C_2}}
\frac{(s^\ddagger)^{\frac{3}{2}} \Xi(s)^{\frac{2}{3}}}{bc \sqrt{C_1C_2mns}} \Big(1+\frac{MN^2S}{nc R}\Big)^{-\frac{1}{2}}(m, bc)\nonumber\\
&\ll X^\varepsilon H^2M^{1/2}NR 
\sum_{\substack{b\le 2M\\ c\le 2N}} 
\sum_{\substack{C_1\\ C_2}} 
\sum_{n\le N_0, \, s\sim S} 
\frac{(s^\ddagger)^{\frac{3}{2}} \Xi(s)^{\frac{2}{3}}}{c^2 \sqrt{C_1C_2n^3b^3s}}\Big(1+\frac{MN^2S}{nc R}\Big)^{-\frac{1}{2}},\label{UB:B22_C}
\end{align}
where we used Lemma \ref{lm:gcd}. Clearly, the function $g(s):=(s^\ddagger)^{\frac{3}{2}} \Xi(s)^{\frac{2}{3}}$ is multiplicative and
$$
g(p^{\nu}) = \begin{cases}
p  & (\nu=1),
\\
p^{\frac{3}{2}}           & (\nu=2),
\\
p^{\frac{2\nu}{3}}  & (\nu\ge 3).
\end{cases}
$$
It is easy to see that there is an absolute constant $C>0$ such that
$$
\sum_{\nu\ge 0} \frac{g(p^{\nu})}{p^{\nu}}
\le 1+\frac{C}{p} .
$$
Thus
\begin{align*}
\sum_{s\le S} \frac{g(s)}{\sqrt{s}}
& \le S^{\frac{1}{2}} \sum_{s\le S} \frac{g(s)}{s}
\le S^{\frac{1}{2}} \prod_{p\le S} \sum_{\nu\ge 0}\frac{g(p^{\nu})}{p^{\nu}}
\ll S^{\frac{1}{2}}X^\varepsilon.
\end{align*}
Inserting this into \eqref{UB:B22_C} and summing trivially over $b, n, C_1, C_2,c$, we can find
\begin{align*}
\cB_1
\ll X^\varepsilon H^2R\ll X^\varepsilon H^2L^{\frac{1}{2}}.
\end{align*}

\section{Estimate for $\cB_2$}
Following the same initial arguments for $\cB_1$ (and after the same change of variables), we have
\begin{align*}
\cB_2
&\ll X^\varepsilon R^{\lambda-\kappa} \sum_{\substack{b\le 2M\\c\le 2N}} 
\sum_{\substack{\tc_1\le 2N/\tilde{c}\\ \tc_2\le 2N/c}} 
\sum_{\substack{h_1\le H\\ h_2\le H}} 
\sum_{\substack{m\sim M/b\\ n\le \min\{2N/c\tc_1,2N/c\tc_2\}\\ s\sim S}}\\
&\ \ \ \ \ \times
\sum_{\substack{\tn_1\sim N/nc\tc_1\\ \tn_2\sim N/nc\tc_2}} 
\frac{((mn\tn_1\tn_2 s)^\sharp)^{\nu-\kappa}}{(bc\tc_1\tc_2)^{\lambda-\kappa}}(mn^2\tn_1\tn_2 sbc^2\tc_1\tc_2)^\kappa(mn^2\tn_1\tn_2,s^\flat)^\kappa.
\end{align*}
Note that the restriction $[\cR_1]$ becomes
\begin{align*}\tag*{$[\cR_{11}]$}
(mn^2\tn_1\tn_2bc^2\tc_1\tc_2)^{1-2\kappa}\Big(\frac{R}{bc\tc_1\tc_2}\Big)^{2\kappa-2\lambda+1}\geqslant \Big(\frac{s^\flat}{(\Delta,s^\flat)}\Big)^{2\kappa},\end{align*}
which holds automatically if one has
\begin{align*}\tag*{$[\cR_{12}]$}
(MN^2)^{1-2\kappa}R^{2\kappa-2\lambda+1}\geqslant(s^\flat)^{2\kappa},\ \ \lambda\geqslant2\kappa,\ \ \lambda-\kappa\geqslant\frac{1}{2}.\end{align*}

In practice, $\nu-\kappa$ usually takes very small values, we thus assume $\nu=\kappa$ to simply the arguments. Hence we may derive that
\begin{align*}
\cB_2
&\ll X^\varepsilon (MN^2S)^\kappa R^{\lambda-\kappa} \sum_{\substack{b\le 2M\\c\le 2N}} 
\sum_{\substack{\tc_1\le 2N/\tilde{c}\\ \tc_2\le 2N/c}} 
\sum_{\substack{h_1\le H\\ h_2\le H}} 
\sum_{\substack{m\sim M/b\\ n\le \min\{2N/c\tc_1,2N/c\tc_2\}\\ s\sim S}}
\sum_{\substack{\tn_1\sim N/nc\tc_1\\ \tn_2\sim N/nc\tc_2}} 
\frac{1}{(bc\tc_1\tc_2)^{\lambda-\kappa}}\\
&\ll X^\varepsilon H^2(MN^2S)^{1+\kappa} R^{\lambda-\kappa}\\
&\ll X^\varepsilon H^2L^{(\lambda+1)/2}(MN^2)^{\kappa+1}.
\end{align*}

To simplify $[\cR_{11}]$, we introduce stronger restrictions
\begin{align*}\tag*{$[\cR_{12}]$}
(MN^2)^{1-2\kappa}R^{2\kappa-2\lambda+1}\geqslant L^{\kappa},\ \ \lambda\geqslant2\kappa,\ \ \lambda-\kappa\geqslant\frac{1}{2}.\end{align*}
Note that $R\asymp L^{1/2},$ then $[\cR_{12}]$ can be reformulated (up to a positive constant) as
\begin{align*}\tag*{$[\cR_{13}]$}
(MN^2)^{1-2\kappa}\geqslant L^{\lambda-\frac{1}{2}},\ \ \lambda\geqslant2\kappa,\ \ \lambda-\kappa\geqslant\frac{1}{2}.\end{align*}

Following the similar manner, we may find that restrictions in $[\cR_2]$
will become 
\begin{align*}\tag*{$[\cR_{21}]$}
(MN^2)^{2-2\lambda}\geqslant L^{\kappa},\ \ \lambda\geqslant2\kappa,\ \ \lambda-\kappa\geqslant\frac{1}{2}\end{align*}
parallel with $[\cR_{13}]$.

\section{Conclusion}

Collecting the above estimates, we find
\begin{align*}
\cB_1&\ll X^\varepsilon (LHMN+H^2L^{\frac{1}{2}})\end{align*}
and
\begin{align*}
\cB_2&\ll X^\varepsilon (L^{(\lambda-\kappa+1)/2}HMN+H^2L^{(\lambda+1)/2}(MN^2)^{\kappa+1})\end{align*}
subject to the restriction $[\cR_{13}]$ or $[\cR_{21}]$. Hence we conclude from \eqref{eq:B(M,N)-decomposition} that
\begin{align*}
\cB(M, N)&\ll H^3N+X^\varepsilon (HLMN+H^2L^{(\lambda+1)/2}(MN^2)^{\kappa+1}).
\end{align*}

\bibliographystyle{plain}

\end{document}